\documentclass[
framed_theorems,
framed_proofs,
print,
]{OCPreprint}

\usepackage{enumitem}

\renewcommand\d{\mathop{}\!\mathrm{d}}

\title[A convex function with empty subdifferential]{A convex, finite and lower semicontinuous function with empty subdifferential}
\author{Gerd Wachsmuth\footnote{%
		Brandenburgische Technische Universität Cottbus--Senftenberg,
		Institute of Mathematics,
		03046 Cottbus,
		Germany,
		\email{wachsmuth@b-tu.de},
		\url{https://www.b-tu.de/fg-optimale-steuerung/team/prof-gerd-wachsmuth},
		ORCID: 0000-0002-3098-1503%
	}
}
{
	\makeatletter
	\def\and{ and }
	\def\footnote#1{}
	\hypersetup{
		pdftitle={\@title},
		pdfauthor={\@author}
	}
	\makeatother
}

\begin{document}
\maketitle
\begin{abstract}
	We give an example of a convex, finite and lower semicontinuous function
	whose subdifferential is everywhere empty.
	This is possible since the function is defined on an incomplete normed space.
	The function serves as a universal counterexample
	to various statements in convex analysis in which completeness is required.
\end{abstract}


\begin{msc}
	\mscLink{46N10},
	\mscLink{90C25}
\end{msc}

\section{Introduction}
\label{sec:intro}
Some results in convex analysis require
that the underlying space is complete and that the (convex) functions
are lower semicontinuous.
As examples,
we mention
\begin{itemize}
	\item
		the Brøndsted--Rockafellar theorem
		about the density of the domain of the subdifferential
		in the domain of the function,
		see
		\cite[Theorem~2]{BrondstedRockafellar1965},
	\item
		maximal monotonicity of the subdifferential,
		see
		\cite[Theorem~A]{Rockafellar1970:2},
	\item
		Ekeland's variational principle,
		see
		\cite[Theorem~1.1]{Ekeland1974},
	\item
		strong Fenchel duality,
		see
		\cite[Theorems~17, 18]{Rockafellar1974}
		and
		\cite[Corollary~1]{Robinson1976},
	\item
		the formula for the convex conjugate of a sum,
		\cite[Theorem~(1.1)]{AttouchBrezis1986}.
\end{itemize}
It is clear that some of these results are closely connected.
The Brøndsted--Rockafellar theorem is usually
proved via Ekeland's variational principle.
Similarly,
strong Fenchel duality
is 
intimately related to the convex conjugate of a sum
and to the sum rule for the subdifferential.

A natural question is whether the assumptions of completeness and lower semicontinuity
in the above results are actually necessary.
If one drops the lower semicontinuity,
linear unbounded functionals often serve as a counterexample.
In absence of completeness,
counterexamples are typically harder to construct,
although it is stated
``It is easy to see that, in general, the conclusion of Theorem (1.1) fails if $E$ is a (non-complete) normed space''
in
\cite{AttouchBrezis1986},
but no concrete example is given.

We are mainly interested in the question of non-emptiness of the subdifferential.
For a normed space $X$ and a proper, convex function $f \colon X \to (-\infty, \infty]$,
the subdifferential of $f$ at $x \in X$
is defined via
\begin{equation*}
	\partial f(x)
	:=
	\set{
		x\dualspace \in X\dualspace
		\given
		\forall y \in X :
		f(y) \ge f(x) + \dual{x\dualspace}{y - x}
	}
	.
\end{equation*}
For many examples, one can check that there always exists some
$x \in X$ such that $\partial f(x) \ne \emptyset$.
In fact,
we already mentioned the
Brøndsted--Rockafellar theorem, which ensures
that
\begin{equation*}
	\dom(\partial f)
	:=
	\set{
		x \in X
		\given
		\partial f(x) \ne \emptyset
	}
\end{equation*}
is dense in
\begin{equation*}
	\dom(f)
	:=
	\set{
		x \in X
		\given
		f(x) < \infty
	}
	,
\end{equation*}
whenever
$X$ is complete and $f$ is lower semicontinuous.
If $f$ is an unbounded linear functional,
it is clear that $\dom(\partial f) = \emptyset$
is not dense in $\dom(f) = X$.
Moreover,
in \cite{BrondstedRockafellar1965}
an example of a proper, convex and lower semicontinuous function $f$
defined on an incomplete space such that $\dom(\partial f) = \emptyset$
is given,
by building upon an example by \cite{Klee1958}.
This construction, however, is quite involved.
We are not aware of similar examples in the literature.

The example in \cite{Rainwater1988},
see also \cite[Example~3.8]{Phelps1993},
comes close, since it possesses an empty subdifferential on a dense subset of $\dom(f)$.
However, this example is posed in the Hilbert space $\ell^2$
and, therefore, the subdifferential cannot be empty at every point
due to the Brøndsted--Rockafellar theorem.

Finally, we mention that
the assumptions of completeness of $X$ and
lower semicontinuity of $f$
can often be replaced by
continuity at a single point of some involved function
and, in this case, completeness of $X$ is not necessary.
A famous example is the validity of the sum rule
\begin{equation*}
	\partial f(x) + \partial g(x)
	=
	\partial (f + g)(x),
\end{equation*}
for convex functions $f,g \colon X \to (-\infty, \infty]$
under the Moreau--Rockafellar condition,
i.e.,
whenever there exists a point $x_0 \in X$
\begin{equation}
	\label{eq:moreau-rockafellar}
	\cont(f) \cap \dom(g) \ne \emptyset,
\end{equation}
where $\cont(f) \subset \dom(f)$ is the set
of all points at which $f$ is continuous.
We emphasize that this result does not require
that $X$ is complete
or that any of the functions $f$ and $g$ is lower semicontinuous.
To complete the picture,
we mention that the sum rule also holds
provided that
\begin{subequations}
	\label{eq:rockafellar-robinson}
	\begin{align}
		\label{eq:rockafellar-robinson:1}
		&
		0 \in \core( \dom(f) - \dom(g) ),
		\\
		\label{eq:rockafellar-robinson:2}
		&
		\text{$X$ is complete and $f, g$ are lower semicontinuous}
	\end{align}
\end{subequations}
are satisfied.
Here, $\core(A)$ is the algebraic interior of a set $A \subset X$.
The condition \eqref{eq:rockafellar-robinson:1}
is called the Rockafellar--Robinson condition
and it is easy to check that it is implied by \eqref{eq:moreau-rockafellar}.
However, \eqref{eq:rockafellar-robinson:2}
is not needed for the sum rule if \eqref{eq:moreau-rockafellar}
is satisfied.
Thus, \eqref{eq:rockafellar-robinson} is not weaker than \eqref{eq:moreau-rockafellar}.
A natural question is whether the sum-rule still holds
provided that only \eqref{eq:rockafellar-robinson:1} is satisfied.
Again, unbounded linear functionals show that the lower semicontinuity assumption
cannot be dropped.
The necessity of the completeness of $X$ is slightly harder to verify.
One possibility is to use a convex
and lower semicontinuous function $f \colon X \to \R$
with an empty subdifferential and check
\begin{equation*}
	X\dualspace
	=
	\partial( f + \delta_{x} )(x)
	\ne
	\partial f(x) + \partial\delta_{x}(x)
	=
	\emptyset + X\dualspace
	=
	\emptyset
\end{equation*}
for some arbitrary $x \in X$.
Due to $\dom(f) = X$,
condition \eqref{eq:rockafellar-robinson:1}
is satisfied.
It remains to construct such a function $f$.

\section{The function and its properties}
\label{sec:function}
Let us consider the linear space of real-valued, finite sequences
\begin{equation*}
	c_c
	:=
	\set{
		\seq{x_n}_{n \in \N} \subset \R
		\given
		\exists N \in \N :
		\forall n > N :
		x_n = 0
	}
	.
\end{equation*}
An easy application of Baire's theorem yields
that $c_c$ equipped with any norm will be an incomplete space.
We equip $c_c$ with the norm of $\ell^2$, i.e.,
\begin{equation*}
	\norm{x}
	:=
	\parens[\bigg]{
		\sum_{n = 1}^\infty x_n^2
	}^{1/2}
	\qquad\forall x \in c_c
	.
\end{equation*}
Note that the dual space of $c_c$
(equipped with this norm)
can be canonically identified with $\ell^2$.
We consider the function
$f \colon c_c \to \R$
defined via
\begin{equation*}
	f(x)
	:=
	\sum_{n = 1}^\infty \frac{n^2}{2} \parens*{ x_n - n^{-2}}^2
	\qquad\forall x \in c_c
	.
\end{equation*}

This function has some very peculiar properties.
\begin{theorem}
	\label{thm:prop}
	Consider the space $c_c$ and the function $f \colon c_c \to \R$ as above.
	\begin{enumerate}[label=(\alph*)]
		\item\label{thm:prop:1}
			The function $f$ is well defined, convex and lower semicontinuous.
		\item\label{thm:prop:2}
			The function $f$ is nowhere continuous.
		\item\label{thm:prop:3}
			At every $x \in c_c$
			the function $f$ is directionally differentiable
			and
			\begin{equation*}
				f'(x; y)
				=
				\sum_{n = 1}^\infty
				n^2 \parens*{x_n - n^{-2}} y_n
			\end{equation*}
			is the directional derivative in direction $y \in c_c$.
			In particular, the
			directional derivative at $x$
			is an unbounded, linear functional.
		\item\label{thm:prop:4}
			At every $x \in c_c$,
			we have $\partial f(x) = \emptyset$.
	\end{enumerate}
\end{theorem}
\begin{proof}
	\ref{thm:prop:1}:
	The function $f$ is well defined,
	since $x_n$ is zero for large $n$
	and since $\sum_{n = 1}^\infty n^{-2}$ converges.
	The convexity is clear
	and the lower semicontinuity
	follows from the lemma of Fatou.

	\ref{thm:prop:2}:
	For every $x \in c_c$ and $n \in \N$ large enough,
	we have
	\begin{equation*}
		f(x + 2 n^{-1} e_n)
		-
		f(x)
		=
		\frac{n^2}{2}
		\parens*{
			\parens*{ 2 n^{-1} - n^{-2} }^2
			-
			\parens*{ - n^{-2} }^2
		}
		\ge
		\frac{n^2}{2}
		\parens*{
			n^{-2}
			-
			n^{-4}
		}
		\ge
		\frac14
		.
	\end{equation*}
	Since $x + 2 n^{-1} e_n \to x$,
	the function $f$ cannot be continuous at
	the arbitrary point $x \in c_c$.

	\ref{thm:prop:3}:
	The formula for the directional derivative
	is clear
	since $x$ and $y$ are finite sequences.
	Consequently,
	$f'(x; \cdot)$ is linear
	and
	$f'(x; e_n) = -1$
	for large enough $n$.
	Thus, $f'(x; \cdot)$ is unbounded.

	\ref{thm:prop:4}:
	From the well-known formula
	\begin{equation*}
		\partial f(x)
		=
		\set{
			z \in \ell^2
			\given
			\forall y \in X :
			\dual{z}{y} \le f'(x; y)
		}
	\end{equation*}
	the emptiness follows,
	since $f'(x; \cdot)$
	is linear and unbounded.
\end{proof}

We emphasize that the verification of the properties of
$f$ is very elementary.

The theorem already shows that
the assertion of the Brøndsted--Rockafellar theorem
is not valid for $f$, thus its completeness assumption is crucial.
Similarly, one can check that the assertion of
the Ekeland variational principle
(which is at the heart of the Brøndsted--Rockafellar theorem)
cannot hold,
since this would give rise to a non-empty subdifferential
of $f$
at some point in $c_c$.
Similarly, the subdifferential $\partial f$ is, of course,
not maximally monotone,
i.e., the completeness of the space is necessary in \cite[Theorem~A]{Rockafellar1970:2}.
Further, the function $f$ cannot be reconstructed from its subdifferential,
since for every linear and continuous functional $\ell$,
we have $\partial f = \partial (f + \ell)$,
but $f$ and $f + \ell$ do not differ by a constant.
Again, completeness in
\cite[Theorem~B]{Rockafellar1970:2}
is crucial.

It is also well known,
see, e.g., \cite[Corollary~8B]{Rockafellar1974},
that convex and lower semicontinuous functions on a Banach space
are continuous in the interior of their domain.
The function $f$ demonstrates that this is not true in incomplete spaces.

Next, we check that the duality results are no longer valid.
Indeed, we just use $g = \delta_{\set{0}}$,
i.e., the indicator function of the origin.
Then, the Fenchel dual of
\begin{equation*}
	\text{Minimize}
	\qquad
	f(x) + g(x)
	\qquad\text{with respect to }
	x \in c_c
\end{equation*}
is
\begin{equation*}
	\text{Maximize}
	\qquad
	- f\conjugate(y) - g\conjugate(-y)
	\qquad\text{with respect to }
	y \in \ell^2.
\end{equation*}
Note that $g\conjugate \equiv 0$.
If the dual problem would possess a solution $y \in \ell^2$,
we would have
\begin{equation*}
	0 \in \partial f\conjugate(y)
\end{equation*}
and, consequently, $y \in \partial f(0)$,
but this is impossible.
Thus, the dual problem does not have a solution.
However,
one can check that no duality gap occurs.

Similarly,
one can check that
\begin{equation*}
	\partial (f + g)(0)
	=
	\ell^2
	\ne
	\emptyset
	=
	\partial f(0) + \partial g(0),
\end{equation*}
i.e., the sum rule fails.
Finally,
\begin{equation*}
	(f + g)\conjugate
	=
	f\conjugate \oplus g\conjugate
	\equiv
	\frac{\pi^2}{12}
\end{equation*}
holds, where ``$\oplus$'' indicates infimal convolution,
but the infimal convolution is not exact,
since the function $f\conjugate$ does not possess a minimizer.

\section{Dual pair with positive duality gap}
\label{sec:gap}
We have seen that a very simple function $g$
is sufficient to get a dual problem without a solution.
We show that using an operator $A \colon c_c \to c_c$
yields a pair of problems with positive duality gap.

We define $A \colon c_c \to c_c$ via
\begin{equation*}
	(A x)_n
	:=
	\begin{cases}
		x_1 & \text{if } n = 1, \\
		x_n - x_{n - 1} & \text{if } n > 1
	\end{cases}
\end{equation*}
for all $x \in c_c$ and $n \in \N$.
It is clear that $A$ is linear
and the boundedness is easy to check,
since
\begin{equation*}
	\norm{A x}^2
	=
	\sum_{n = 1}^\infty (A x_n)^2
	\le
	x_1^2
	+
	\sum_{n = 2}^\infty \parens*{ 2 x_n^2 + 2 x_{n-1}^2}
	\le
	4 \norm{x}^2
	.
\end{equation*}
Thus $\norm{A} \le 2$
and one can check that we actually have $\norm{A} = 2$.

Further, let
$g \colon c_c \to \R$ be the zero function.
We start by investigating the primal problem
\begin{equation*}
	\text{Minimize} \quad
	f(A x) + g(x)
	\qquad\text{with respect to }
	x \in c_c
	.
\end{equation*}
We check that $0$ is the solution of this problem.
In fact,
for every $y \in c_c$,
the directional derivative of $f \circ A$ at $0$
in direction $y$
equals
\begin{equation*}
	(f \circ A)'(0; y)
	=
	f'(0; A y)
	=
	-y_1
	+
	\sum_{n = 2}^\infty 
	- ( y_n - y_{n - 1} )
	=
	0
	,
\end{equation*}
see \itemref{thm:prop:3}.
Since $f$ is convex and $g \equiv 0$, it follows that $0$
is a solution of the primal problem
and the primal infimal value is $\pi^2/12$.

It is also interesting to note that this implies
\begin{equation*}
	\partial( f \circ A )(0)
	=
	\set{0}
	\ne
	A\adjoint \partial f( A 0 )
	=
	\emptyset
	,
\end{equation*}
i.e., the failure of the chain rule.

Now, we consider the dual problem
\begin{equation*}
	\text{Maximize}
	\qquad
	- f\conjugate(y) - g\conjugate(-A\adjoint y)
	\qquad\text{with respect to }
	y \in \ell^2.
\end{equation*}
Note that $g\conjugate = \delta_{\set{0}}$
and a short calculation shows
\begin{equation*}
	f\conjugate(y)
	=
	\sum_{n = 1}^\infty \frac{1}{n^2} \parens*{\frac12 y^2 + y}
	,
\end{equation*}
since we can argue coefficient-wise.
Since the operator $A$ has a dense range,
the operator $A\adjoint$ is injective
(which can also be shown directly).
Consequently,
$y = 0$ is the only point with
$g\conjugate(-A\adjoint y) < \infty$.
Thus, this is the solution of the dual problem
and the objective value is $0$.
Due to
\begin{equation*}
	0
	<
	\frac{\pi^2}{12}
	,
\end{equation*}
this pair of problems possesses a positive duality gap.

We mention that
$\dom f = c_c$ implies
\begin{equation*}
	0 \in \core(A \dom g - \dom f ).
\end{equation*}
Further,
the functions $f$ and $g$ are lower semicontinuous.
Only the non-completeness of $c_c$
causes problems.
This prevents \cite[Theorem~18]{Rockafellar1974}
from being applied.

There also does not exist a point $x_0 \in \dom g$
such that $f$ is continuous at $A x_0$,
since $f$ is not continuous at $0$ and $\dom g = \set{0}$.
Consequently, \cite[(8.25)]{Rockafellar1974}
is not satisfied.

\section{Conclusion}
Using elementary arguments,
we have verified
in \cref{thm:prop}
that the subdifferential of the function $f$ is everywhere empty.
Consequently,
this function shows that the completeness assumption in many theorems of convex analysis
cannot be dropped.


\printbibliography

\end{document}